\documentclass[onecolumn]{autart}    % Enable this line and disable the 
\usepackage{graphicx}
\usepackage[utf8]{inputenc}
\usepackage{mathptmx} % assumes new font selection scheme installed
\usepackage{newtxtext}
\usepackage{textcomp}
\usepackage{enumitem}
\usepackage{dsfont, url}
\usepackage{mathtools}
\usepackage[usenames, dvipsnames]{xcolor}
\usepackage[%
colorlinks={true},%
citecolor={Blue},%
urlcolor={PineGreen},%
linkcolor={Sepia}%
]{hyperref}
\usepackage{float}
\usepackage{newfloat}
\usepackage[skip=1pt]{caption}
\usepackage{amsmath} % assumes amsmath package installed
\usepackage{amssymb}
\usepackage{adjustbox}
\usepackage{verbatim} 
\usepackage{xargs}
\usepackage{cite}

\DeclareMathOperator*{\minimize}{minimize}
\DeclareMathOperator*{\sbjto}{subject\ to}
\DeclareMathOperator*{\blkdiag}{bdiag}

\DeclareMathOperator{\rank}{rank}

\DeclareMathOperator{\trace}{tr}

\renewcommand{\leq}{\leqslant}

\renewcommand{\geq}{\geqslant}

\newcommand{\R}{\mathds{R}}
\newcommand{\Nz}{\mathds{N}_0}
\newcommand{\N}{\mathds{Z}_{+}}
\newcommand{\EE}{\mathds{E}}

\newcommand{\bmat}[1]{\begin{bmatrix}#1\end{bmatrix}}
\newcommand{\abs}[1]{\left|#1\right|}
\newcommand{\norm}[1]{\left\|#1\right\|}
\newcommand{\secref}[1]{\S \ref{#1}}

\newcommand{\transp}{^\top}

\newcommand{\zeros}{\mathbf{0}}

\newcommand{\st}{x}

\newcommand{\stest}{\tilde{x}}

\newcommand{\dortho}{d_o}

\newcommand{\A}{A}
\newcommand{\calA}{\mathcal{A}}
\newcommand{\Aortho}{\A_o}
\newcommand{\Aschur}{\A_s}
\newcommand{\B}{B}

\newcommand{\Bortho}{\B_o}
\newcommand{\Bschur}{\B_s}

\newcommand{\control}{u}

\newcommand{\wnoise}{w}

\newcommand{\cnoise}{\nu}
\newcommand{\snoise}{s}

\newcommand{\calB}{\mathcal{B}}
\newcommand{\calD}{\mathcal{D}}

\newcommand{\calG}{\mathcal{G}}
\newcommand{\calH}{\mathcal{H}}
\newcommand{\calQ}{\mathcal{Q}}
\newcommand{\calR}{\mathcal{R}}

\newcommand{\calS}{\mathcal{S}}

\newcommand{\XX}{\mathfrak X}

\newcommand{\ee}{\mathfrak{\psi}}
\newcommand{\offset}{\boldsymbol{\eta}}
\newcommand{\gain}{\boldsymbol{\Theta}}
\newcommand{\authority}{u_{\max}}
\newcommand{\reachindex}{\kappa}
\newcommand{\reachab}{\mathrm{R}}
\newcommand{\Let}{\coloneqq}
\newcommand{\teL}{\eqqcolon}

\allowdisplaybreaks
\definecolor{mycolor1}{rgb}{0.00000,0.44700,0.74100}%
\definecolor{mycolor2}{rgb}{0.85000,0.32500,0.09800}%

\begin{document}

\begin{frontmatter}

\title{Reference tracking stochastic model predictive control over unreliable channels and bounded control actions} % Title, preferably not more 

\author[UIUC]{Prabhat K. Mishra}\ead{pmishra@illinois.edu}, 
\author[MIT]{Sanket S. Diwale}\ead{ diwales@mit.edu},
\author[EPFL]{Colin N. Jones}\ead{colin.jones@epfl.ch},
\author[IIT]{Debasish Chatterjee}\ead{dchatter@iitb.ac.in},

\address[UIUC]{Coordinated Science Laboratory, University of Illinois at Urbana-Champaign (UIUC), USA.}
\address[MIT]{Massachussets Institute of Technology (MIT), USA.}
\address[EPFL]{Automatic Control Laboratory, École Polytechnique Fédérale de Lausanne (EPFL), Switzerland. }
\address[IIT]{Systems \& Control Engineering,
	Indian Institute of Technology Bombay (IITB), 
	India.}        

\begin{abstract}                          % Abstract of not more than 200 words.
A stochastic model predictive control framework over unreliable Bernoulli communication channels, in the presence of unbounded process noise and under bounded control inputs, is presented for tracking a reference signal. The data losses in the control channel are compensated by a carefully designed transmission protocol, and that of the sensor channel by a dropout compensator. A class of saturated, disturbance feedback policies is proposed for control in the presence of noisy dropout compensation. A reference governor is employed to generate trackable reference trajectories and stability constraints are employed to ensure mean-square boundedness of the reference tracking error. The overall approach yields a computationally tractable quadratic program, which can be iteratively solved online.  
\end{abstract}

\end{frontmatter}
\section{Introduction}
\par Regulation and tracking are two basic problems in control theory and have been studied rigorously in the framework of model predictive control (MPC) \cite{ref:rawlings-09}. A tracking problem can be considered as the generalization of a regulation problem with a time-varying set point. However, with a varying set point, the stabilizing and recursively feasible design of a regulating MPC may not remain valid \cite{limon2013Encyclopedia}. Consequently, the  tracking problem generally requires a separate treatment. 
\par The stability and recursive feasibility issues for tracking MPC in linear, disturbance-free systems have been tackled in works like \cite{Borrelli2016TAC,Limon2008tracking,periodictracking_jones} using constructions of appropriate control invariant sets around a reachable pseudo trajectory. Since networked control systems are gaining prominence due to their flexibility, the reference tracking problem for networked systems without additive process noise is considered in \cite{pang2014outputtracking}. The above references \cite{Borrelli2016TAC, Limon2008tracking, periodictracking_jones, pang2014outputtracking} give a good understanding of the underlying challenges and techniques to solve tracking problems but do not consider process noise.% in 
\par Tracking MPC under bounded process noise is considered in \cite{limon2010robusttracking, paulson2018mixed, output_tracking}. In many applications, a hard bound on the additive process noise is not known a priori, but a distribution with unbounded support can be safely prescribed. In such cases, a \emph{stochastic MPC} (SMPC) formulation is desired. \cite{Mesbah2018Tracking, farina2017tracking} are two examples of such SMPC formulations with input chance constraints. However, for physical systems, with hard input constraints, providing only chance constraints guarantees for input bounding may cause input saturation and instability in the physical system. A SMPC formulation for tracking with hard constraints on control actions is as yet missing in the literature.
\par It is well known for regulation problems in SMPC that a class of feedback policies is essential for control over the uncertainty \cite{kumar1986stochastic, mayne2016_right_direction, goulart-06}. In order to satisfy hard constraints on control actions, \cite{ref:ChaHokLyg-11} employed saturated disturbance feedback policies and ensured stability by using stability constraints for a sufficiently large control authority. Another such class of policies with evolving saturated disturbances is employed in \cite{Tsiotras_covariance_hard} for covariance steering. The approach in \cite{ref:ChaHokLyg-11} was extended for any positive control authority in \cite{ref:PDQ-15} and the controller was implemented on networks with the help of novel transmission protocols. 
\par It is argued in \cite{ref:PDQ-15} that a stabilizing SMPC formulation over networks requires three ingredients -- a class of feedback policies,  appropriate transmission protocols and stability constraints. This framework is utilized in our previous contributions \cite{ref:PDQ-15, PDQ_Policy} assuming an unreliable control channel (uplink) and a perfect sensor channel (downlink). However, with an unreliable downlink, the computation of disturbance feedback of the above works is made impossible and cannot be applied for such situations. We fill this particular lacuna in this article. 
\par Generally, in the case of incomplete and corrupt measurements, a Kalman filter is used for estimating the missing signals with some degree of uncertainty. The saturated innovation term of the Kalman filter can then be used in the feedback policy as in \cite{ref:Hokayem-12, PDQ-LCSS}. This approach was used under the settings of an unreliable downlink in \cite{PDQ_intermittent} by using a Kalman filter on the plant side and computing feedback from the Kalman filter innovation terms, sent through the downlink. This, however, does not fully address the problem of communication dropouts in the downlink and requires addition computational resources on the plant side for the filtering. Furthermore, all the approaches above \cite{ref:ChaHokLyg-11, ref:PDQ-15, ref:Hokayem-12, PDQ_intermittent, PDQ-LCSS} present a constrained SMPC approach for \emph{a regulation problem} and their adaptation under the setting of a time varying reference signal is also missing in the literature.
\par In this article, we adapt and extend the results of \cite{ref:PDQ-15} by considering an unreliable downlink and a time-varying set point, reference tracking problem. A class of feedback policies in which the feedback term takes the additive uncertainties appearing in the dynamics of the dropout compensator is considered. We show that the proposed class of feedback policies leads to a tractable  constrained optimal control problem and present modifications to the transmission protocols and stability constraints presented in \cite{ref:PDQ-15} under the settings of bounded-input reference tracking and unreliable communication in order to get a provably stabilizing and tractable SMPC framework. 
\par The contributions of this article are (a) SMPC formulation of tracking with hard constraints on control actions and (b) a new class of feedback policies in the context of an unreliable downlink that ensures tractability of the underlying constrained stochastic optimal control problem. 
\par This article proceeds as follows: We present the problem formulation and system setup in \secref{s: problem setup}. The main ingredients of the tracking problem under unreliable channels are explained in \secref{s: reference governor}, \secref{s: dropout compensator} and \secref{s: controller}, respectively. We present our main result on tractability and stability in \secref{s: tractability and stability}. We illustrate our results in \secref{s: numerical experiment} by numerical experiments and conclude in \secref{s: epilogue}.
\subsection*{Notation}
We use the symbol \( \zeros \) to denote a matrix of appropriate dimensions with all elements $0$ and \(I_d\) for the \(d\times d\) identity matrix. We simply use $I$ in place of \(I_d\) when its dimensions are clear from the context. The notation \(\EE_z[\cdot]\) is used for the conditional expectation given \(z\). Let $\sigma_1(M)$ denote the largest singular value of $M$, and $M^{\dagger}$ its Moore-Penrose pseudo inverse. A block diagonal matrix $M$ with diagonal entries $M_1, \cdots , M_n$ is represented as $M = \blkdiag\{M_1, \cdots , M_n\}$. For a vector $V \in \R^d$, $V^{(i)}$ denotes its $i^{th}$ entry for $i = 1, \ldots, d$. For any vector sequence \((v_n)_{n\in\Nz}\), let \( v_{n:k} \) denote the vector \(\bmat{v_n\transp & v_{n+1}\transp & \cdots & v_{n+k-1}\transp}\transp\), \(k\in\N\). 

\section{Problem setup and solution architecture}\label{s: problem setup}

%\subsection{Plant}
We consider a discrete-time dynamical system
\begin{equation}\label{e:system}
\st_{t+1} = \A \st_t + \B \control_t^a + \wnoise_t,
\end{equation}
where $\st_t \in \R^d$, $\control_t^a \in \R^m$, $\wnoise_t \in \R^d$, are the system state, input and process noise, respectively. The matrix pair $(\A, \B)$ is assumed to be controllable. The additive process noise \((\wnoise_t)_{t\in\Nz}\) is assumed to be a sequence of i.i.d. zero mean random vectors taking values in \(\R^d\). Each component of $\wnoise_t$ is symmetrically distributed about the origin and $\wnoise_t$ has bounded fourth moments, i.e., $\EE [\norm{\wnoise_t}^4] < \infty$.\footnote{The bounded fourth moment assumption being much weaker than a Gaussian distribution assumption, allows for more general noise distribution characteristics in this work.}
\par
\begin{figure}
	\centering
	\begin{adjustbox}{width=0.8\columnwidth}
		\includegraphics{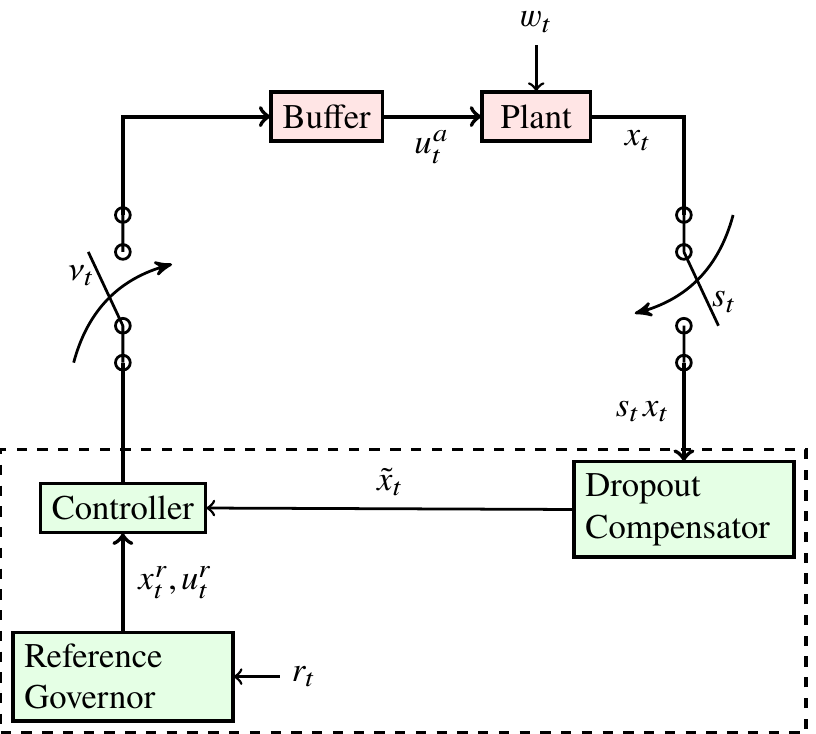}		
	\end{adjustbox}	
	\caption{A dropout compensator feeds the estimated states to the controller and the controller generates the admissible control sequence by taking reference governor into account. The computed control parameters are transmitted through the control channel by a suitably chosen transmission protocol and the successfully received parameters are stored in a buffer.}
	\label{fig:blockDia}
\end{figure}

 The control input $\control_t^a$ is required to be bounded and without loss of generality is assumed to be uniformly bounded as
\begin{equation}\label{e:controlset}
\norm{\control_t^a}_{\infty} \leq \authority \text{ for all } t \in \Nz.
\end{equation} 
We use the superscript notation $\control_t^a$ in the system dynamics to denote the actual applied control action to the system which may be different from the computed control input $u_t$ due to an unreliable communication channel (uplink); See Fig.\ \ref{fig:blockDia}. 
\subsection{Communication channels}
\par The controller communicates with the dynamical system (plant) through unreliable Bernoulli channels. The parameters of the control policy are transmitted through a control channel (uplink) from the controller to the plant. State information $x_t$ is transmitted through a sensor channel (downlink) from the plant to the controller. Successful transmissions from an unreliable channel follow a Bernoulli distribution. To model this distribution, we utilize two i.i.d. Bernoulli random variable sequences $(\nu_t)_{t \in \Nz}$ with mean $p_c$ and $(s_t)_{t \in \Nz}$ with mean $p_s$ for the control channel and the sensor channel, respectively. A transmission across the control channel and the sensor channel is considered successful if $\nu_t=1$ and $s_t=1$, respectively. Each successful transmission from controller to the actuator is accompanied by an acknowledgment of success to the sender as done in TCP like protocols.  
\subsection{Controller} 
\par The aim of a tracking controller in the presence of perfect channels is to design a control sequence $(\control_t)_{t \in \Nz}$ such that the state process $(\st_t)_{t\in \Nz}$ converges to the given reference sequence $(r_t)_{t \in \Nz}$. Since $(r_t)_{t \in \Nz}$ may not be trackable due to the dynamical constraints, generally it is approximated by a trackable signal $(\st_t^r)_{t \in \Nz}$ typically obtained by a reference governor. Section \ref{s: reference governor} details the design of the reference governor.

Even in the presence of perfect channels the tracking problem is not trivial when the process noise has unbounded support and control inputs are uniformly bounded. Since the traditional form of asymptotic stability cannot be achieved in the presence of additive stochastic uncertainties, we address the notion of mean-square boundedness. We recall the following definition:
\begin{defn}\cite[\S III.A]{chatterjee-15}\label{def:msb}
	{\rm
		An $\R^d$-valued random process $(\st_t)_{t \in \Nz}$ is said to be mean square bounded with respect to the available information $\XX_0$ at $t=0$ at the controller if there exists $\gamma < \infty$ such that $ \sup_{t \in \Nz}\EE_{\XX_0}[\norm{\st_t}^2] \leq \gamma $. 
	}
\end{defn}  
The above notion of mean-square boundedness implies that the probability of the norm of $x_t$ being greater than or equal to $k$, decays faster than $k^{-2}$, uniformly over time \cite[III.A]{chatterjee-15}. Moreover, this notion also implies that there does not exist any divergent trajectory, with probability one \cite[Lemma 4.1]{thesis_milan_masters}. Therefore, for the purpose of tracking in a stochastic environment, we focus on designing a controller such that the underlying error between the state of the system and the reference signal remains mean-square bounded.
\subsection{Solution architecture}
\par Given a reference signal $(r_t)_{t \in \Nz}$, the reference governor in \secref{s: reference governor}, generates sequences $(\st_t^r)_{t \in \Nz}$ and $(\control_t^r)_{t \in \Nz}$ such that the deviation of $x_t^r$ from $r_t$, denoted as the \emph{governor error} $e_t^{G} = \st_t^r - r_t$ is bounded. The reference control input $\control_t^r$ is made to satisfy a bound $\norm{\control_t^r}_\infty\leq \delta u_{max}$ for some $\delta \in(0,1)$ to limit the control authority available for tracking and noise compensation respectively; see \secref{s: reference governor} for details.
\par The reference governor signals $\st^r_t,\control_t^r$ are utilized by the SMPC controller to generate an open-loop control sequence $\offset_t$ and disturbance feedback gain matrix $\gain_t$; such that the closed-loop reference tracking error, denoted $e_t=\st_t-\st_t^r$ is guaranteed to be mean square bounded. However, due to the unreliable sensor channel, the state vector $\st_t$ is not always available for feedback computation in the controller. Instead, a state estimator is used as a dropout compensator and a state estimate $\stest_t$ is made available to the controller. The tracking error $e_t$, can thus be split into two errors, (i) the estimator error $e_t^D=\st_t - \stest_t$ and (ii) the controller error $e_t^C=\stest_t - x_t^r$.
The dropout compensator is designed in \secref{s: dropout compensator} such that the estimation error $e_t^{D} = \st_t - \stest_t$ remains mean-square bounded. Further, the SMPC controller is designed in \secref{s: controller} such that the controller error $e_t^{C} = \stest_t - \st_t^r$ remains mean-square bounded. Thus, the \emph{overall error} between the state of the system $x_t$ and the reference signal $r_t$ written as, 
\begin{equation} \label{e:overall error}
e_t^O = e_t + e_t^G = e_t^{D} + e_t^{C} + e_t^{G}
\end{equation} 
remains mean square bounded, due to each component error being mean-square bounded (Theorem \ref{th:stability}). 
We conclude by showing that the overall approach is computationally tractable and satisfies hard constraints on control actions in \secref{s: tractability and stability}.  The error signals are summarized in Table \ref{T:error summary}.

\begin{minipage}{\linewidth}
	\centering
	\captionof{table}{Error signals} \label{T:error summary} 
	\begin{tabular}{c |  c | c} 
		\hline
		error & expression & governing equations \\  
		\hline
		tracking error & $e_t = \st_t - \st_t^r$  & \eqref{e:tracking_error_redefined}\\ 
		governor error & $ e_t^{G} =\st_t^r - r_t$ &\eqref{e:bound reference governor} and \eqref{e:governor_error}\\
		estimator error & $e_t^{D} = \st_t - \stest_t$ & \eqref{e: DC error}\\
		controller error & $e_t^{C} = \stest_t - \st_t^r$  & \eqref{e: controller error}\\ 
		overall error & $e_t^O = \st_t - r_t$     & \eqref{e:overall error}\\
		[1ex] 
		\hline
	\end{tabular}
\end{minipage}

The following sections proceed with further details of the control system architecture components, viz, the reference governor, dropout compensator and SMPC controller.

\section{Reference governor}\label{s: reference governor}
For a given constrained linear system only a class of signals may be trackable \cite{gorecki2015, periodictracking_jones}. This class of trackable signals is generally obtained by a reference governor \cite{borrelli2009RG, stoican2012reference}. Although tracking of arbitrary reference signals is not possible with bounded control inputs, a  trackable pseudo reference trajectory is feasible with the help of a reference governor. The tracking performance can then be monitored with respect to this generated reference. In subsequent analysis, we focus on a class of trackable signals as defined below:
\begin{defn}\label{e:rGovernor}
	\rm{
		For a given control authority $\authority$, weight $\delta \in (0,1)$ and error bound $\gamma^{G}\geq 0$, a signal $(r_t)_{t \in \Nz}$ is called \emph{trackable} if there exist sequences $(\st_t^r)_{t \in \Nz}, (\control_t^r)_{t \in \Nz}$ such that
		\begin{equation} \label{e:bound reference governor}
		\left.
		\begin{aligned}
		& \st_{t+1}^r = \A \st_t^r + \B \control_t^r, \\
		& \sup_{t \in \Nz} \left( \norm{\st_t^r - r_t}^2 \right) \leq \gamma^{G}, \\
		& \norm{\control_t^r}_{\infty} \leq \delta \authority \text{ for each } t.  
		\end{aligned}
		\right\}
		\end{equation}	
	}
\end{defn}
The weight $\delta$ is used to distribute the control authority between the governor  input $\control_t^r$ and SMPC controller for additional disturbance rejection.The difference between the trackable trajectory $\st_t^r$ and reference signal $r_t$ is denoted as the \emph{governor error} 
\begin{equation}\label{e:governor_error}
e^G_t = \st_t^r - r_t.
\end{equation}
The sequences $(\st_t^r)_{t \in \Nz}, (\control_t^r)_{t \in \Nz}$ can be obtained by solving an optimal control problem of the form
\begin{equation} \label{e: reference governor OCP}
\begin{aligned}
(\st_t^r)_{t\leq T}, (\control_t^r)_{t\leq T} \; =&\quad \underset{(\st_t),(\control_t)}{\text{argmin}} \sum_{t=0}^T \norm{ \st_t - r_t }^2 + \norm{\control_t}^2\\
\sbjto &\quad    \st_{t+1} = \A \st_t + \B \control_t\\
& \quad \norm{\control_t}_{\infty} \leq \delta \authority; \qquad \forall \; t\leq T\\
& \quad \st_0 = r_0.
\end{aligned}	
\end{equation}
Given $r_t$ and $\delta$, $r_t$ is considered admissible if it provides a feasible solution to \eqref{e: reference governor OCP} and a finite bound $\gamma^G$, to render $\st_t^r$ trackable according to Definition \ref{e:rGovernor}. The problem \eqref{e: reference governor OCP} can then be solved offline for a finite horizon problem or online using a moving horizon, for infinite horizon reference signals as done in \cite{Farina2013mpc_rg}. We refer the reader to \cite{tutorial_reference_governor, Farina2013mpc_rg, Bemporad98governor} for other advances in reference governor designs. 

Given a trackable pair of sequences, $(\st_t^r,\control_t^r)$, our objective boils down to ensuring mean square boundedness of the \emph{tracking error}, $e_t$, 
by using the control authority $(1-\delta)\authority$ and present a tractable formulation of the SMPC controller. 

\section{Dropout compensator}\label{s: dropout compensator}
At each time step, the state information $\st_t$ is transmitted from the plant to the controller over an unreliable communication channel as described in Section \ref{s: problem setup}. A Bernoulli random variable $s_t$ determines the successful transmission of data across the channel. For the event, $s_t=1$ (with probability $p_s$), $\st_t$ is successfully transmitted and for $s_t=0$ (with probability $1-p_s$), the transmission fails.

When a successful transmission is received, $s_t=1$, the state estimate for the controller is set to the received value, i.e., $\stest_t = \st_t$. For failed transmission, $s_t=0$, the previous state estimate is propagated forward using the system dynamics
\begin{equation}
\stest_{t}=A\stest_{t-1} + B u_{t-1}^a
\end{equation}
initialized with $\stest_{-1}=\zeros$ and $u^a_{-1}=\zeros$. For $t\geq0$, $u_t^a$ is given by the disturbance feedback policy \eqref{e:policy} discussed in \secref{s: controller}.
The combined state estimate update can then be written as,
\begin{equation}\label{e:est}
\stest_t = \snoise_t  \st_t + (1-\snoise_t) \left( \A \stest_{t-1} + \B \control_{t-1}^a \right).
\end{equation}

Dropout compensators of the form \eqref{e:est} are widely used in the literature \cite{gupta2007optimal} and a justification for \eqref{e:est} can be found in \cite[Lemma 3]{PDQ_intermittent}. Using \eqref{e:system} and \eqref{e:est}, the estimator error dynamics can then be written as,
\begin{equation}\label{e: DC error}
e_t^{D} \Let \st_t - \stest_t = (1-\snoise_t)(\A e_{t-1}^{D} + \wnoise_{t-1}).
\end{equation}

\section{Controller}\label{s: controller}
Let $N, N_r \in \N$ be two positive integer constants such that $N_r \leq N$. The stochastic MPC controller solves a $N$-step stochastic optimal control problem after every $N_r$ steps. Recalling that the controller objective is to minimize the reference tracking error,
\begin{equation}\label{e:tracking_error_redefined}
e_t \Let \st_t - \st_t^r,
\end{equation}
the objective function for the SMPC controller is set to
\begin{equation}\label{e:cost}
V_t \Let \EE_{\XX_{t}} \left[e_{t+N} \transp Q_f e_{t+N} +  \sum_{i = 0}^{N-1} \left( e_{t+i} \transp Q e_{t+i} + (\control_{t+i}^e) \transp R \control_{t+i}^e \right) \right],
\end{equation}  
where $\control_{t+i}^e \Let \control_{t+i}^a - \control_{t+i}^r $ and $\XX_t$ is the information available at the controller up to and including time $t$. Matrices $Q$, $Q_f$ are taken to be symmetric, positive semi-definite and $R$ is taken to be positive definite. The controller transmits the parameters of the policy (defined in \secref{s:policy}) over the control channel with the help of a transmission protocol (defined in \secref{s:tp}) such that previously transmitted parameters can be used in the case of transmission failures.
The backbone of SMPC is to iteratively solve the following constrainted stochastic optimal control problem (CSOCP) in the presence of stability constraints (defined in \secref{s:stability constraint}) over the class of policies: 
%\begin{equation}\label{e:csocp}
%\left .
%\begin{aligned}
%\underset{\text{policy}}{\minimize} & \quad V_t \\
%\sbjto & \quad \eqref{e:system}, \eqref{e:controlset}, \eqref{e:drift intermittent}.
%\end{aligned}
%\right \}
%\end{equation}
\begin{equation}\label{e:csocp}
\underset{\text{policy}}{\minimize}  \quad \{ V_t \mid \eqref{e:system}, \eqref{e:controlset}, \eqref{e:drift intermittent} \}.
\end{equation}
The above CSOCP \eqref{e:csocp} is re-solved after every $N_r$ time steps with a new initial state estimate and the updated decision variables are transmitted through the control channel to allow the computation of $u^a_t$ at the plant input by using the transmission protocol. We discuss these mechanisms in the following subsections.
\subsection{Class of feedback policies}
\label{s:policy}
The presence of buffer and unreliable control channel is generally ignored at the time of controller design. In this article, we follow the approach of \cite{ref:PDQ-15} to incorporate both effects at the design stage itself. For this purpose, we have different notations and equations for the control policy and the applied control to the system, $\control_t$ and $\control_t^a$, respectively.
The dynamics of the dropout compensator \eqref{e:est} is rewritten as
\begin{equation}\label{e: compensator dynamics}
\stest_{t+1} = \A \stest_t + \B \control_t^a + \tilde{\wnoise}_t,
\end{equation}
with $\tilde{\wnoise_t} \Let \snoise_{t+1}\left( \A e_t^{D} + \wnoise_t \right)$.
The class of causal policies is then parametrized in terms of the compensator disturbances $\tilde{\wnoise_t} $ as:
\begin{equation}\label{e:policy_temp}
\control_{t+i} = \control^r_{t+i} + \eta_{t+i} + \sum_{j=0}^{i} \theta_{i,j}\ee(\tilde{\wnoise}_{t+j -1})
\end{equation}
where $\control^r_{t+i}$ is the reference control input computed in \secref{s: reference governor}; $(\eta_{t+i})_{i={0,\dots,N-1}}$ and $(\theta_{i,j})_{j\leq i: i=0,\dots,N-1}$ are the nominal input and feedback gain respectively, to be  computed by solving CSOCP \eqref{e:csocp}. The estimator disturbance $\tilde{\wnoise}_{t+j -1}$ is computed from \eqref{e: compensator dynamics} as
\begin{equation}\label{e:compensator_noise}
\tilde{\wnoise}_{t+j -1} = \stest_{t+j} - (A\stest_{t+j-1} + B\control^a_{t+j-1}).
\end{equation}
The computation is initialized with $\stest_{-1}=\zeros$, $\control^a_{-1}=\zeros$ (as per \secref{s:  dropout compensator}).
The saturation function $\ee:\R \rightarrow \R$, can be any anti-symmetric function; such that, $\ee(0)=0$, $\ee(-x)=-\ee(x)$ and $\sup_{x\in\R}\ee(x) \leq \ee_{\max}$. For vector inputs, like $\tilde{\wnoise}_{t+j -1} $, $\ee$ is applied element-wise.
The class of policies \eqref{e:policy_temp} can be written in a compact matrix form as:
\begin{equation}\label{e:policy}
\control_{t:N} = \control_{t:N}^r + \offset_t + \gain_t \ee(\tilde{\wnoise}_{t-1:N}),
\end{equation} 
where $\offset_t \Let \eta_{t:N}$, \(\offset_t\in\R^{m N}\) and \(\gain_t\) is a lower block triangular matrix
\begin{equation} \label{e:gain}
\gain_t = \bmat{ \theta_{0, t} & \zeros & \cdots & \zeros & \zeros \\ \theta_{1, t} & \theta_{1, t+1}  & \cdots & \zeros & \zeros \\ \vdots & \vdots & \vdots & \vdots & \vdots \\ \theta_{N-1, t} & \theta_{N-1, t+1} & \cdots & \theta_{N-1, t+N-2} & \theta_{N-1, t+N-1} },
\end{equation}
with each \(\theta_{k, \ell} \in \R^{m\times d}\) and $\norm{\ee(\tilde{\wnoise}_{t:N-1})}_{\infty} \leq \varphi_{\max}$.
\subsection{Transmission protocol}\label{s:tp}
Recalling that $\control_{t:N}^r$ is obtained from the reference governor, and that $(\offset_t, \gain_t)$ is the solution of the CSOCP \eqref{e:csocp} solved at time $t$, $\control_{t:N}^r, \offset_t$ and $\gain_t$ are available at time $t$. These can be used to compute control inputs $\control_{t:N}$ with the help of the causally available information of $\tilde{w}_t$. The absence of the future feedback terms at the optimization instant, motivates the selection of the transmission rule of the policy parameters as in \cite{ref:PDQ-15}. 
The idea of transmission protocols is inspired by the so-called packetized predictive control techniques \cite{Quevedo-12}, where a buffer stores a finite sequence of future control actions to be applied in case of control channel failures. We initialize a buffer to hold $N_r$ number of control inputs and adapt a transmission protocol of \cite{ref:PDQ-15} below:
\begin{enumerate}[label= {(TP\arabic*)}, leftmargin= *, widest=3, align=left, start=1, nosep]
	\item \label{a:repetitive} At each time $t\geq 0$, do:
	\begin{enumerate}[ leftmargin=*, widest=3, align=left, start=1]
		\item If $t=kN_r$, for some $k\in \Nz$,
		\begin{enumerate}
			\item Compute $\control_{t:N}^r$, $\offset_{t}$ and $\gain_{t}$ using the optimization based controller from Section \ref{s: controller}
			\item Compute $\control_{t:N_r}$ using \eqref{e:policy}
			\item Set $\ell=0$
			\item Transmit $\left\{\control_{t} , (\control_{t+1:N_r-1}^r+\eta_{t+1:N_r-1})\right\}$ to the buffer
		\end{enumerate}  
		\item else,
		\begin{enumerate}
			\item Update $\ell=\ell+1$
			\item If buffer is empty, transmit\\
			$ \left\{\control_{t+\ell}, (\control_{t+\ell+1:N_r-\ell-1}^r+\eta_{t+\ell+1:N_r-\ell-1}) \right\} $ \\
			Otherwise, transmit $\control_{t+\ell}$.
		\end{enumerate}
	\end{enumerate}	
\end{enumerate}	
As mentioned in \ref{a:repetitive}, $N_r-\ell$ blocks of control are transmitted at each time $t+\ell$, when buffer is empty. Otherwise only one block of the control is transmitted. If the transmitted block is received at the buffer, then it is stored there starting from the beginning. At each time instant the first block of the buffer is applied to the plant and then the buffer is rotated by using a left shift register. This operation of shift register makes the buffer automatically empty before each optimization instant. 
\par Let $g_t = \cnoise_t$, $g_{t+\ell} = g_{t+\ell-1} + (1- g_{t+\ell-1})\cnoise_{t+\ell}$. The term $g_{t+\ell}$ then gives state of consecutive dropouts in the unreliable control channel and
the transmission protocol \ref{a:repetitive}, thus, effectively applies an input to the actuator given by
\begin{equation}
\control_{t+\ell}^a = g_{t+\ell}(\control_{t+\ell}^r + \eta_{t+\ell}) + \cnoise_{t+\ell}\sum_{i=0}^{\ell} \theta_{\ell,t+i}\ee (\tilde{\wnoise}_{t+ i -1})
\end{equation}	
at time $t+\ell$. 
Let us define a block diagonal matrix $\calS \Let \blkdiag\{\cnoise_{t}I_{m}, \ldots, \cnoise_{t + \reachindex -1}I_{m}, I_{m(N-\reachindex)}\}$ and the matrix $\calG$, which has $N \times N$ blocks in total, each of dimension $m \times m$ and for $i,j = 1,\cdots,N$, the matrix $\calG$ is given in terms of the blocks $\calG_b^{(i,j)}$ each of dimension $m \times m$ as 		 
\begin{align}
\calG_b^{(i,j)} \Let
\begin{cases}
g_{t+i-1} I_m & \text{ if } i=j \leq N_r ,\\ 
I_{m} & \text{ if } i =j > N_r ,\\
0_m \quad  & \text{ otherwise.}			
\end{cases}		
\end{align}
The stacked control vector in one optimization horizon can then be rewritten as:
\begin{equation}\label{e:policyrepetitive out}
\control^a_{t:N} \Let \calH \control_{t:N}^r + \calG \offset_t + \calS \gain_t  \ee(\tilde{\wnoise}_{t-1:N}), \\
\end{equation}
where $\gain_t$ and $\ee(\wnoise_{t:N+1})$ are as defined in \eqref{e:policy} and $\calH = \calG$. 
\begin{rem}
	\rm{
		In any application where the future reference control signals $\control_t^r$, which are generated by the reference governor, are known a priori, we can transmit the future values of $\control_t^r$ a priori and store in a buffer on the plant side. In such situations we can assume that $\control_t^r$ is available at the actuator at time $t$. Therefore, we can set $\calH = I$. In the present article, we do not consider this apriori storage, which sets $\calH = \calG$.
	}
\end{rem} 	
\subsection{Stability constraints}\label{s:stability constraint}
Similar to \cite{ref:HokChaRamChaLyg-10}, stability of the proposed approach is independent of the cost function. We employ stability constraints, which are derived from \cite{ref:ChaRamHokLyg-12} under the following assumptions:
\begin{enumerate}[label= {\rm(A\arabic*)}, leftmargin= *, widest=3, align=left,  nosep]	
	\item \label{as:processnoise} The zero mean noise sequence $(\wnoise_t)_{t \in \Nz}$ is fourth moment bounded, i.e., $\EE [\norm{\wnoise_t}^4] \leq C_4$, for some $C_4 < \infty$.
	\item \label{as:Lyapunov} The system matrix $\A$ has all eigenvalues in the unit disk and those on unit circle are semi-simple\footnote{The algebraic and geometric multiplicities are same for the eigenvalues that are on the unit circle.}.	
	\item \label{as:controlability} System matrix pair $(A,B)$ is controllable.
\end{enumerate}		
The above assumptions \ref{as:processnoise} - \ref{as:Lyapunov} are also used in \cite{ref:HokChaRamChaLyg-10} along with stabilizability of the matrix pair $(A,B)$ for the case of perfect channels \footnote{\ref{as:Lyapunov} refers to the largest class of linear time invariant systems known to be stabilizable by bounded control actions under perfect communication channels. Please refer to \cite[Abstract]{Sussmann97} and \cite[Theorem 1.7]{ref:ChaRamHokLyg-12} for more details.}. We need \ref{as:controlability} because mere stabilizability is not sufficient for reference tracking problems. We quickly recall some steps so that the stability constraints along the lines of \cite{PDQ_intermittent} can be employed under the settings of the present article. Let us first define the \emph{controller error} $e_t^{C} \Let \stest_t - \st_t^r$, then from \eqref{e: compensator dynamics} and \eqref{e:bound reference governor} we get
\begin{equation}\label{e: controller error}
e_t^{C} \Let \stest_t - \st_t^r = \A e_{t-1}^{C} + \B \control_{t-1}^e + \tilde{\wnoise}_{t-1},
\end{equation}
where $ \control_{t}^e = \control_{t}^a - \control_{t}^r$.
Without loss of generality due to the Assumption \ref{as:Lyapunov}, we can assume that the error dynamics \eqref{e: controller error} is of the form
\begin{equation} \label{e:orthogonal decomposition error}
\bmat{(e^{C}_{t+1})^o \\ (e^{C}_{t+1})^s} = \bmat{\Aortho & 0\\ 0 & \Aschur}\bmat{(e^{C}_{t})^o \\ (e^{C}_{t})^s}	+ \bmat{\Bortho\\ \Bschur}\control_t^e + \bmat{\tilde{\wnoise}_t^o \\ \tilde{\wnoise}_t^s}, 
\end{equation}
where \(\Aortho\in\R^{d_o\times d_o}\) is orthogonal and \(\Aschur\in\R^{d_s\times d_s}\) is Schur stable, with \(d = d_o + d_s\). By the controllability assumption \ref{as:controlability}, there exists a positive integer $\reachindex$ such that the reachability matrix \[ \reachab_{\reachindex}(\Aortho, \Bortho) \Let \bmat{\Aortho^{\reachindex - 1}\Bortho & \cdots & \Aortho\Bortho & \Bortho} \] has full row rank; i.e., $\rank(\reachab_{\reachindex}(\Aortho, \Bortho)) = \dortho$. We present our tracking problem in such a way that \cite[Lemma 7]{PDQ_intermittent} is applicable in the context of this article and provides the following result: 
\begin{lem}\label{t:msbsingle}
	Given a system of the form \eqref{e:orthogonal decomposition error} and assumptions \ref{as:processnoise}-\ref{as:Lyapunov}, if for $t=0,\reachindex, 2\reachindex, \ldots$, the following conditions hold: 
	\begin{subequations} \label{e:drift intermittent}
		\begin{align}
		&	\Bigl( (\Aortho^{t+ \reachindex })\transp\reachab_{\reachindex}(\Aortho, \Bortho)\EE_{\XX_t} [\control_{t:\reachindex}^e] \Bigr)^{(j)} \leq -\zeta \notag  \\
		& \quad \text{ whenever }  \left( (\Aortho^t) \transp (e_t^C)^o \right)^{(j)} > c, \label{e:drift1 intermittent} \\	
		&	\Bigl( (\Aortho^{t+  \reachindex})\transp\reachab_{\reachindex}(\Aortho, \Bortho)\EE_{\XX_t} [\control_{t:\reachindex}^e]\Bigr)^{(j)} \geq \zeta \notag \\
		& 	\quad  \text{ whenever }  \left( (\Aortho^t) \transp (e_t^C)^o \right)^{(j)} < -c, \label{e:drift2 intermittent}
		\end{align}
	\end{subequations}
	for each $j = 1, \cdots, \dortho$, for some $\zeta, c>0$;	
	then under the transmission protocol \ref{a:repetitive} there exists $\gamma^C> 0$ such that 
	\begin{equation}\label{e:bound on controller error}
	\EE_{\XX_0} \left[ \norm{e^{C}_t}^2 \right] \leq \gamma^C \quad \text{ for all } t.
	\end{equation} 
	Moreover, when $\zeta \in \left]0, \frac{(1-\delta)\authority}{\sqrt{d_o}\sigma_{1}\left(\reachab_{\reachindex}(\Aortho,\Bortho)^{\dagger}\right)} \right] $, there exists a $\reachindex$-history dependent class of inputs $\control_{\reachindex t : \reachindex}^e$ such that \eqref{e:drift intermittent} and \eqref{e:controlset} are satisfied. 
\end{lem}	
\begin{pf}
	Applying \cite[Lemma 7]{PDQ_intermittent} on \eqref{e:orthogonal decomposition error}, we get \eqref{e:drift intermittent} and the bound \eqref{e:bound on controller error} when $\norm{\control_t^e}_{\infty} \leq (1-\delta) \authority$ by the choice of $\zeta$. The reference governor is designed in \eqref{e:bound reference governor} such that $\norm{\control_t^r}_{\infty} \leq \delta \authority$. Therefore, noting $\norm{\control_t^a}_{\infty} \leq \norm{\control_t}_{\infty} \leq \norm{\control_t^r}_{\infty} + \norm{\control_t^e}_{\infty} \leq \authority$, $\control_t^a$ satisfies \eqref{e:controlset}. 
\end{pf}

\section{Tractability and Stability} \label{s: tractability and stability}
In this section we present how the CSOCP \eqref{e:csocp} can be written as a computationally tractable quadratic program. Recalling $\control_t^e = \control_t^a - \control_t^r$, we get
\begin{equation}\label{e:stacked_control_error}
\control^e_{t:N} \Let (\calH-I) \control_{t:N}^r + \calG \offset_t + \calS \gain_t  \ee(\tilde{\wnoise}_{t-1:N}).
\end{equation}
The compact form representation of error in \eqref{e:cost} over one optimization horizon can be written as:
\begin{equation}\label{e:stacked error}
e_{t:N+1} = \calA e_t + \calB\control_{t:N}^e + \calD \wnoise_{t:N}, 
\end{equation}
where $\calA$, $\calB$, $\calD$ are standard matrices derived by recursing the dynamics over one horizon and given by $\calA = \bmat{I& \A\transp & \cdots & (\A^N)\transp}\transp$,
\[
%\calA = \bmat{I\\ \A\\ \vdots \\ \A^N}, 
\calB = \bmat{\zeros & \zeros & \cdots & \zeros\\ \B & \zeros & \cdots & \zeros \\ \vdots & \vdots & \ddots & \vdots\\ \A^{N-1}\B & \A^{N-2}\B & \dots & \B},  
\calD = \bmat{\zeros & \zeros & \cdots & \zeros\\ I & \zeros & \cdots & \zeros\\ \vdots & \vdots & \ddots & \vdots\\ \A^{N-1} & \A^{N-2} & \dots & I}.
\]
Let $\calQ  \Let \blkdiag\bigl(\overset{N\text{ times}}{\overbrace{Q, \cdots, Q}}, Q_f\bigr), \calR \Let \blkdiag(\overset{N \text{ times}}{\overbrace{R, \cdots, R}})$, the cost function \eqref{e:cost} can also be written in a compact form as follows:
\begin{equation}\label{e:cost_compact}
V_t \Let \EE_{\XX_t} \left[ \norm{e_{t:N+1}}^2_{\calQ} + \norm{\control_{t:N}^e}^2_{\calR} \right].
\end{equation}	 
\par Let us define $\alpha \Let \calB\transp\calQ\calB + \calR$ and the covariance matrices $\Sigma_{\calG} \Let \EE \left[\calG\transp \alpha\calG \right]$, $\Sigma_{\calS} \Let \EE \left[\calS\transp \alpha\calS \right]$,  $\Sigma_{\calH\calG} \Let \EE \left[(\calH-I)\transp \alpha\calG \right]$ and $\Sigma_{\calH\calS} \Let \EE\left[ (\calH-I)\transp \alpha \calS\right]$, where $\calH$ and $\calS$ are defined in \secref{s:tp}. We define two mean matrices $\mu_{\calG} \Let \EE[{\calG}]$ and $\mu_{\calS} \Let \EE[{\calS}] $. The above matrices are needed to write a tractable surrogate of the cost function $V_t$. For a given distribution of $\snoise_t$ and $\cnoise_t$, we can compute them offline by using Monte-Carlo simulations. We need to compute three more covariance matrices $\Sigma_{\ee} \Let \EE_{\XX_t}\Bigl[\ee(\tilde{\wnoise}_{t:N-1})\ee(\tilde{\wnoise}_{t:N-1})\transp \Bigr]$, $\Sigma_{\ee\wnoise} \Let \EE_{\XX_t} \Bigl[\ee(\tilde{\wnoise}_{t:N-1})\wnoise_{t:N} \transp \Bigr]$ and  $\Sigma_{e\ee} = \EE_{\XX_t} \Bigl[ \ee (\tilde{\wnoise}_{t:N-1})(e_t^D) \transp \Bigr]$.
\par Please notice that $\Sigma_{\ee}, \Sigma_{\ee\wnoise}$ and $\Sigma_{e\ee}$ depend on the number of consecutive dropouts of the sensor channel. If we know that there are at most $h$ number of consecutive packet dropouts in the sensor channel, we can compute $h$ number of values of $\Sigma_{\ee}, \Sigma_{\ee\wnoise}$ and $\Sigma_{e\ee}$ offline and use the required one at the time of optimization. Since the sensor channel is assumed to be Bernoulli, the number of consecutive packet dropouts of the sensor channel is unbounded. However, for practical purposes, when data is not received for a long time from a particular channel, we can either give high priority to that channel so that precomputed values of $\Sigma_{\ee}, \Sigma_{\ee\wnoise}$ can be utilized or we compute new values online in advance. The assumption of an uniform bound on the number of consecutive packet dropouts is standard in literature \cite[Assumption 4]{Quevedo2011}. However, we do not need such assumption for any theoretical result in this article.
\par Let
$\gain_t^{(:,t)} \Let \bmat{\theta_{0,t} \transp & \theta_{1,t} \transp & \ldots & \theta_{N-1,t} \transp} \transp $ be the first $d$ columns of $\gain_t$ and $\gain_t^{\prime}$ be such that $\gain_t = \bmat{\gain_t^{(:,t)} & \gain_t^{\prime}}$. Let $\Pi_{\wnoise} = \ee(\tilde{\wnoise}_{t-1}) \ee(\tilde{\wnoise}_{t-1}) \transp $,  We have the following Lemma:
\begin{lem}\label{lem:cost}
	The objective function \eqref{e:cost} can be written as the following convex quadratic function:
	\begin{align}\label{e:obj out channel}
	& V_t^{\prime} =	\offset_t\transp \Sigma_{\calG} \offset_t + 2 \offset_t \transp \Sigma_{\calG\calS}\gain_t^{(:,t)}\ee(\tilde{\wnoise}_{t-1}) 
	+ \trace \Sigma_{\calS} \gain_t^{(:,t)} \Pi_{\wnoise} (\gain_t^{(:,t)})\transp \notag \\ & \quad + \trace (\Sigma_{\calS} \gain_t^{\prime} \Sigma_{\ee} (\gain_t^{\prime})\transp) + 2\trace \left( \calD \transp \calQ \calB\mu_{\calS}\gain_t^{\prime} \Sigma_{\ee\wnoise} \right) \notag \\ & \quad  + 2 \trace \bigl( \calA \transp \calQ \calB \mu_{\calS} \gain_t^{\prime} \Sigma_{e\ee} \bigr) + 2 (e_t^C) \transp \calA \transp \calQ \calB \mu_{\calG}\offset_t \notag  \\
	& \quad + 2(\control_{t:N}^r) \transp \Sigma_{\calH\calS} \gain_t^{(:,t)} \psi(\tilde{\wnoise}_{t-1}) + 2(\control_{t:N}^r) \transp \Sigma_{\calH\calG}\offset_t.
	\end{align} 
\end{lem}
A proof of the Lemma \ref{lem:cost} is given in the appendix. Please notice in the proof of Lemma \ref{lem:cost} that $V_t^{\prime}$ is obtained by removing some constant terms of $V_t$ in \eqref{e:cost}. Therefore, they are not equal but they serve the same purpose for the optimization purpose. Further, the hard constraint on control actions can be written as affine function of the decision variables. We have the following Lemma:
\begin{lem}
	For the class of control policies \eqref{e:policy}, the input constraint \eqref{e:controlset} is equivalent to:
	\begin{equation}\label{e:decisionboundsingle}
	\abs{(\control_{t:N}^r)^{(i)} + \offset_t^{(i)}} +  \norm{\gain_t^{(i,:)}}_1\varphi_{\max}  \leq \authority, 
	\end{equation}	
	for $i = 1, \ldots , Nm$.
\end{lem}	
\begin{pf}
	Since the saturation function is component-wise symmetric about the origin, the claim follows from \cite[Proposition 3]{hokayem2009stochastic} by observing that $\norm{\control_{t:N}^a}_{\infty} \leq \authority \iff \norm{\control_{t:N}}_{\infty} \leq \authority$. 
\end{pf}
The CSOCP \eqref{e:csocp} can be written as the following convex quadratic program:
\begin{equation}\label{e:optimization program}
\minimize \quad  \{ \eqref{e:obj out channel} \mid \eqref{e:decisionboundsingle},~ \eqref{e:drift intermittent} \}.
\end{equation}
The above optimization program \eqref{e:optimization program} can be solved, for example, using the MATLAB based software package YALMIP \cite{ref:Lofberg-04} and the solver GUROBI \cite{gurobi} or SDPT3 \cite{ref:toh-06}. The above optimization program is solved iteratively online after each $N_r$ time steps. The covariance matrices $\Sigma_{\calG}, \Sigma_{\calG\calS}, \Sigma_{\calS}, \Sigma_{\ee}, \Sigma_{\calH\calS}, \Sigma_{\ee\wnoise}, \Sigma_{e\ee}$ and mean matrices $\mu_\calG, \mu_\calS$ are considered as constants in the optimization program \eqref{e:optimization program}. At each optimization instant the estimator computes $\stest_t$ by \eqref{e:est} and $\tilde{\wnoise}_{t-1}$ by \eqref{e:compensator_noise}. The reference signals $\st_t^r$ given by the reference governor is used to compute $e_t^C$ by using \eqref{e: controller error}. The reference control signal $\control_{t:N}^r$ is also obtained from the reference governor. In this way, at each optimization instant $\tilde{\wnoise}_{t-1}, \Pi_\wnoise, e_t^C$ and $\control_{t:N}^r$ are given to the solver to obtain $\offset_t$ and $\gain_t$. We have the following result on mean-square boundedness of the overall error. Its proof is given in the appendix.
\begin{thm}\label{th:stability}
	Consider a discrete-time dynamical system \eqref{e:system} and overall error \eqref{e:overall error} and the control sequence is generated by repeatedly solving the optimization program \eqref{e:optimization program}. Let assumptions \ref{as:processnoise} - \ref{as:controlability} hold, then overall error $e_t^O$ is mean-square bounded.  
\end{thm}	
\section{Numerical Experiments}\label{s: numerical experiment}
\begin{figure*}
	\centering
	\begin{adjustbox}{width = \textwidth}
			\includegraphics{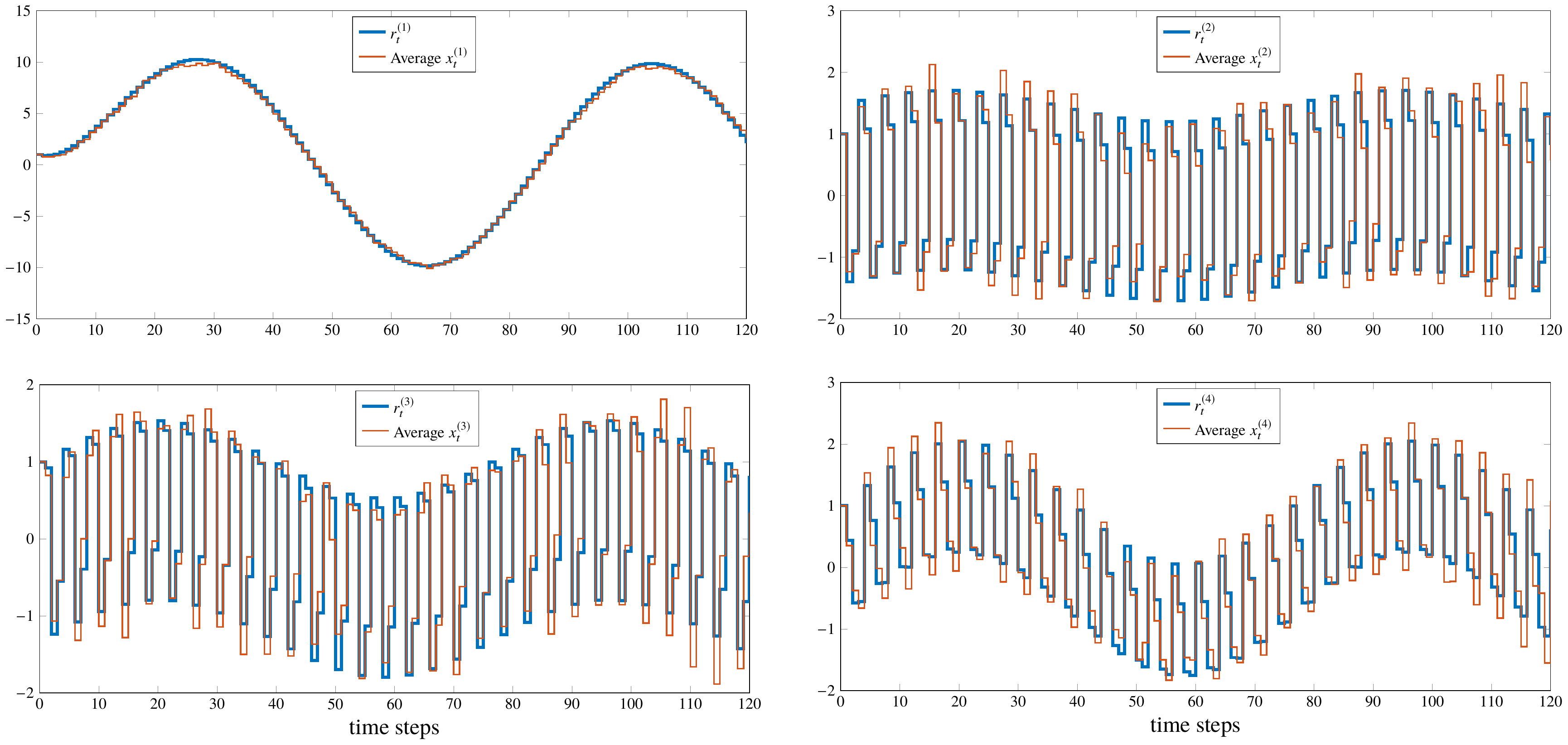}
	\end{adjustbox}
	\caption{States of the reference signal $r_t$ and the averaged states of the system $\st_t$}
	\label{fig:tracking}
\end{figure*}
\begin{figure}
	\centering
	\begin{adjustbox}{width = \columnwidth}
			\includegraphics{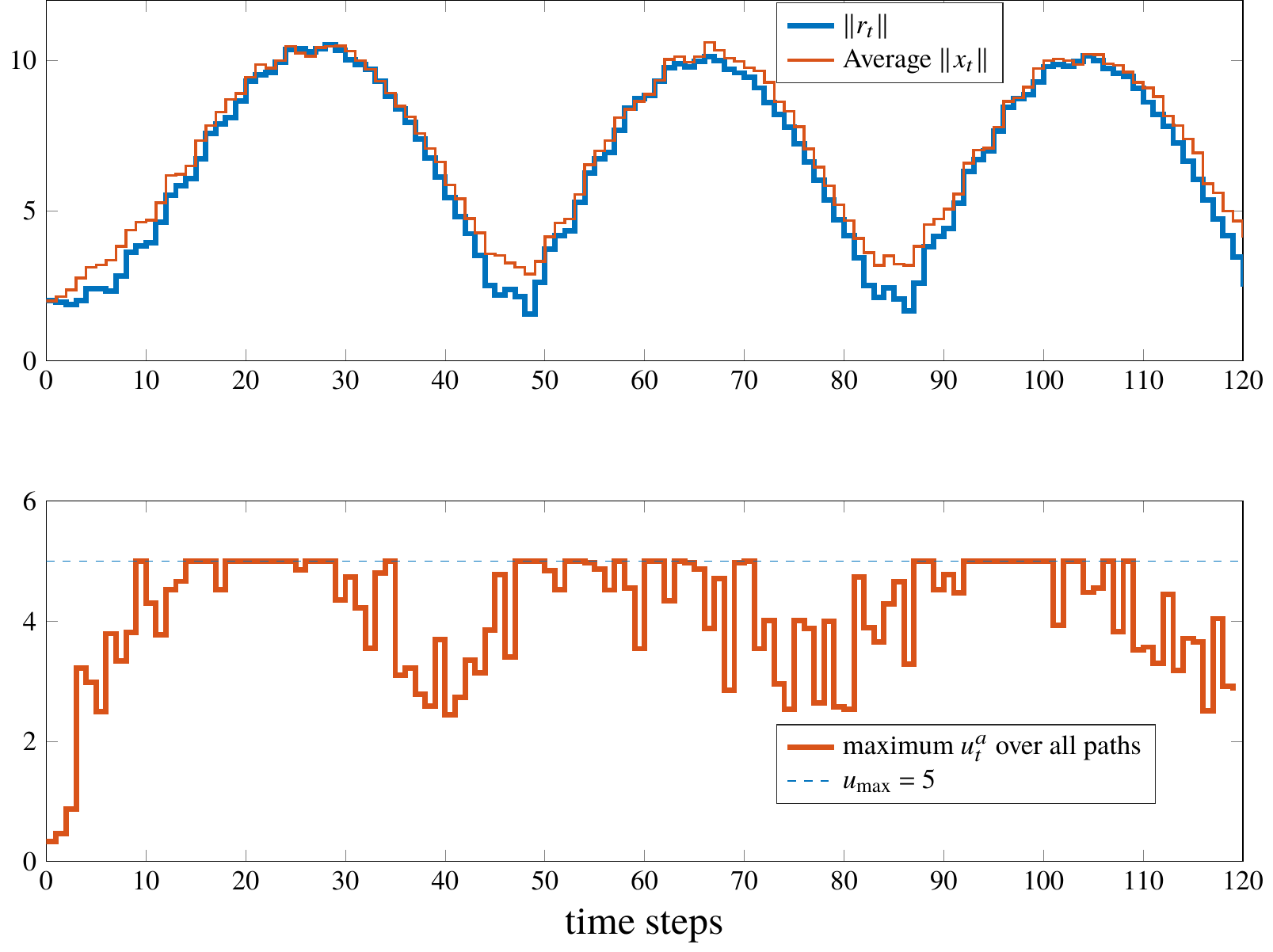}
	\end{adjustbox}
	\caption{Average $\norm{\st_t}$ and maximum control over $50$ sample paths}
	\label{fig:tracking_norm}
\end{figure}
In this section, we present a numerical experiment to validate our theoretical results. We consider a four dimensional stochastic LTI system with system matrices  
\begin{equation*}
\A = \bmat{0.9 & 0 & 0 & 0 \\ 0 & 0 &-0.8 & -0.6\\ 0 & 0.8 & -0.36 & 0.48\\ 0 & 0.6 & 0.48 & -0.64}, 
\B = \bmat{0.5\\ 0.5\\ 0\\ 0.5}, \abs{\control_t^a} \leq 5. 
\end{equation*}
%The control is constrained by the hard bound $\authority = 5$.
The additive process noise is mean zero Gaussian with variance $\Sigma_{\wnoise} = 0.5I$. The successful transmission probabilities for uplink and downlink are 
$p_c = p_s = 0.9$, 
and simulation data is $Q = Q_f= I_4, R =1, N=5, \st_0 = \bmat{1 &1 & 1 & 1}\transp$. The reachability index $\kappa$ of the matrix pair $(\Aortho,\Bortho)$ is $3$. The saturation function  $\ee(\xi) =\frac{1- e^{-\xi}}{1+ e^{-\xi}}$ in our policy is a scalar sigmoidal function which is applied element-wise. We have chosen the recalculation interval $N_r$ same as the reachability index $\reachindex$ of the matrix pair $(\Aortho,\Bortho)$. For simplicity in the reference governor design, we consider a trackable reference signal, which admits the following recursion:
\begin{equation}
r_{t+1} = A r_t + B \left( \frac{5}{2}\sin (0.083 t)\right); \quad r_0 = \st_0 . 
\end{equation}
We choose $\delta = 0.5$, $\st_0^r = r_0$, $\control_t^r = \frac{5}{2}\sin (0.083 t)$, so $\st_t^r = r_t$. For more general reference signals an optimization problem of the form \eqref{e: reference governor OCP} can be solved to compute $\st_t^r, \control_t^r$ or some other methods of reference governor design as in \cite{tutorial_reference_governor, Farina2013mpc_rg, Bemporad98governor} can be adopted. 
\par A simulation for $50$ sample paths of $120$ time steps is performed and the average over all sample paths is demonstrated in Fig.\ \ref{fig:tracking} and Fig.\ \ref{fig:tracking_norm}. We observe in Fig.\ \ref{fig:tracking} that the averaged states follow very closely to the states of the reference signal. Fig.\ \ref{fig:tracking_norm} shows that $\norm{r_t}$ and the averaged $\norm{\st_t}$ are very close while successfully maintaining $\abs{\control_t^a}\leq 5$ in each sample path. 
\par
\begin{figure}
	\centering
	\begin{adjustbox}{width = \columnwidth}
		\includegraphics{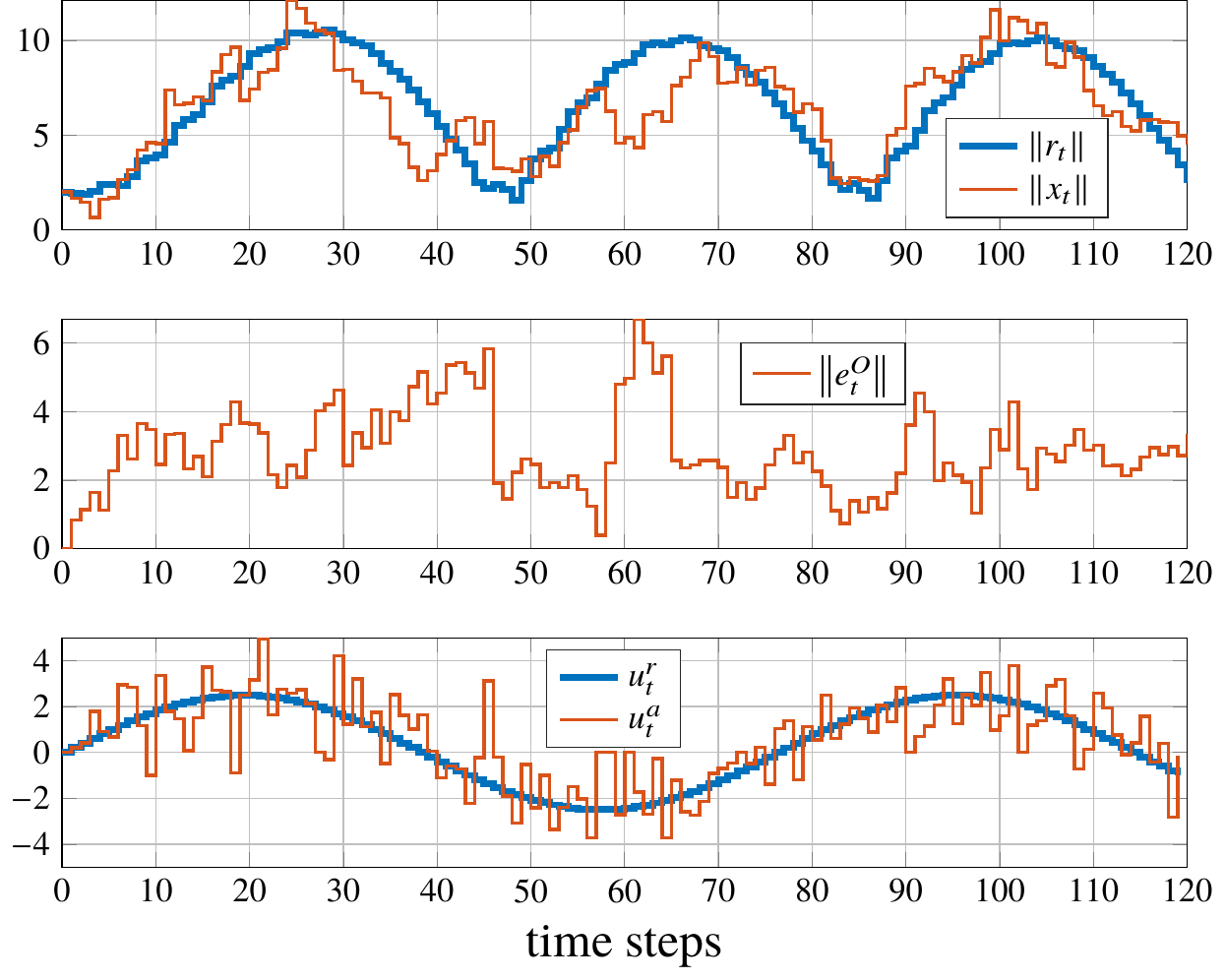}
	\end{adjustbox}
	\caption{Norm of the states, overall error, and control in one path}
	\label{fig:one_state}
\end{figure}

\begin{figure}
	\centering
	\begin{adjustbox}{width = \columnwidth}
		\includegraphics{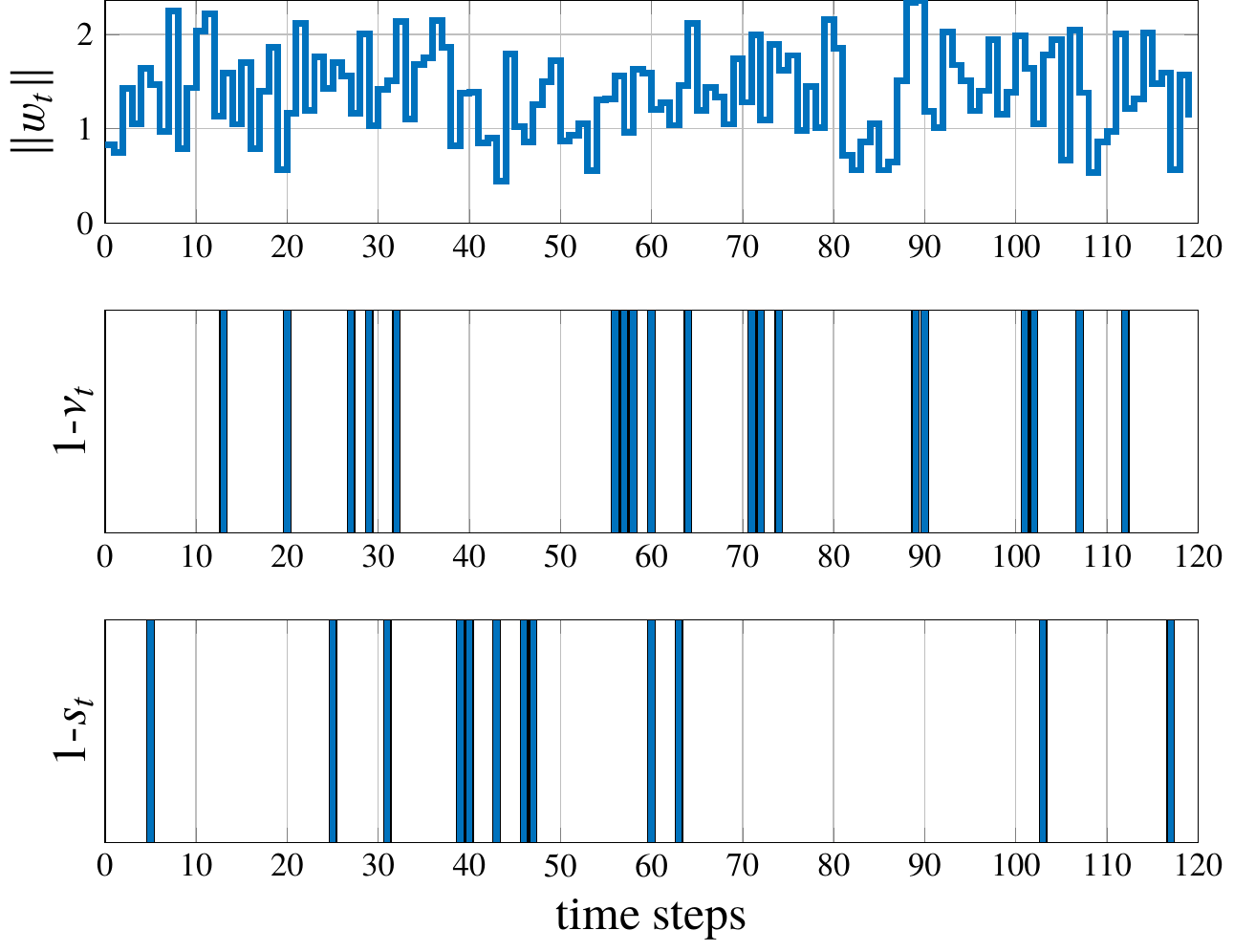}
	\end{adjustbox}
	\caption{Noise in one sample path.}
	\label{fig:noise}
\end{figure}

 Figures \ref{fig:one_state} and \ref{fig:noise} show the data associated with a typical sample path. We observe in Fig.\ \ref{fig:one_state} that $\norm{\st_t}$ and $\norm{r_t}$ closely follow one another in the presence of the additive process noise and packet dropouts in the sensor and control channels, shown in Fig.\ \ref{fig:noise}. The applied control $\control_t^a$ follows the reference control $\control_t^r$ generated by the reference governor with some deviations to compensate for the noise and dropouts. The norm of the overall error $\norm{e_t^O}$ is affected by the additive process noise (see Fig.\ \ref{fig:noise}), and increases when there are consecutive dropouts in one of the channels and is observed to be the highest at $t=61$ when consecutive dropouts of control channel are followed by a dropout at the sensor channel. A sudden decrease in $\norm{e_t^O}$ is observed at $t=46$ due to the decrease in $\norm{\wnoise_t}$, and perfect transmission through sensor channel at $t=45$ and through control channel at $t=46$. Two consecutive dropouts in the sensor channel at $t=46,47$ result in slight increase in $\norm{e_t^O}$ at $t=48$.    
\par
\begin{figure}
	\centering
	\begin{adjustbox}{width = 0.6\columnwidth}
		\includegraphics{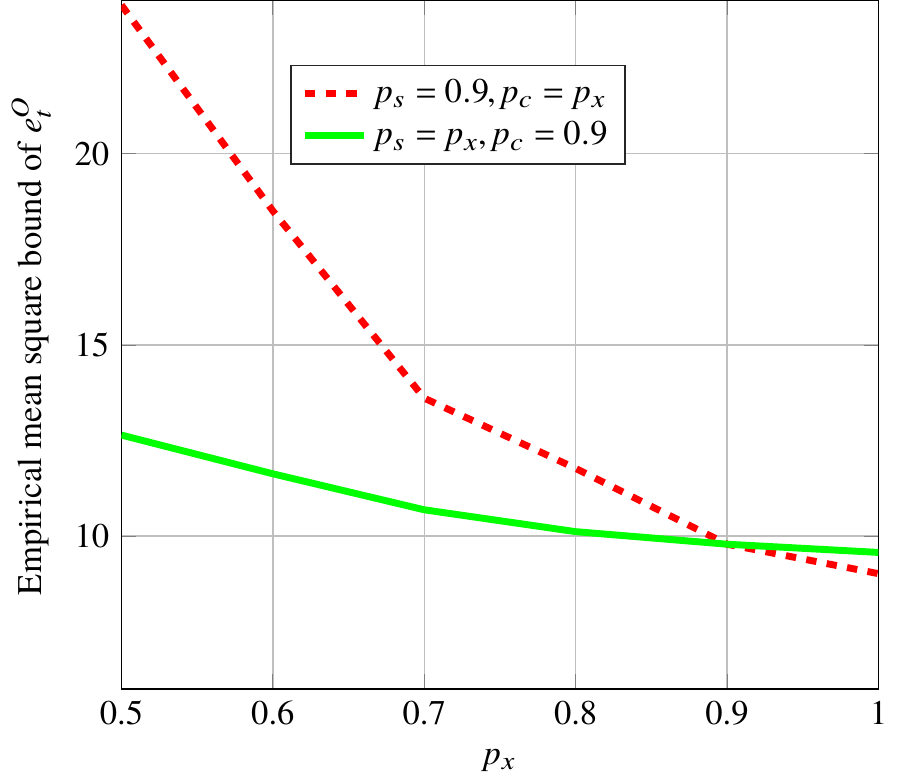}		
	\end{adjustbox}	
	\caption{Empirical mean square bound when dropout probability of one channel varies from 0.5 to 1 while that of the other channel is fixed at 0.9.}
	\label{fig:msb_tracking}
\end{figure}

 Additionally, a parametric study of the tracking performance over wide range of dropout probabilities $p_s$ and $p_c$ is performed. The \emph{empirical mean square bound} (MSB) of the overall error (the error between the states of the controlled plant and the reference signal $e_t^O = \st_t - r_t$) is shown in Fig. \ref{fig:msb_tracking} for varying values of $p_s$ and $p_c$. We first fix $p_c = 0.9$ and record empirical MSB over 200 sample paths when $p_s$ varies in the set $\{0.5, 0.6, 0.7, 0.8, 0.9, 1 \}$. In the second experiment, we fix $p_s = 0.9$ and record the empirical MSB over 200 sample paths when $p_c$ varies in the set $\{0.5, 0.6, 0.7, 0.8, 0.9, 1 \}$. The empirical MSB is larger for fixed $p_s$ when there are more packet dropouts in the control channel. For a fixed $p_c$ also empirical MSB decreases with the increase in $p_s$ but its slope is smaller than that with varying $p_c$ and fixed $p_s$.  


\section{Epilogue}\label{s: epilogue}
The stochastic MPC formulation under unreliable sensor and control channels is considered in this article and the proposed approach is presented for general tracking problems. A class of trackable signals is restricted with the help of a reference governor and available control authority. The work can further be extended in the setting of incomplete and corrupt measurements, and output feedback settings with the help of a Kalman filter along the lines of \cite{PDQ_intermittent}.

\begin{ack}                               % Place acknowledgements
The work of C.\ N.\ Jones has received support from the Swiss National Science Foundation under the Risk project (Risk aware Data Driven Demand Response, grant number 200021 175627).
\end{ack}

\bibliographystyle{plain}        % Include this if you use bibtex 
\bibliography{refs}           % and a bib file to produce the 
                                 % bibliography (preferred). The
                                 % correct style is generated by
                                 % Elsevier at the time of printing.
\appendix
\section{Proofs}    % Each appendix must have a short title.
\begin{pf}[Lemma \ref{lem:cost}]
	Since $e_{t:N+1}$ and $\control_{t:N}^e$ are affine functions of decision variables and the expectation of a convex function is convex \cite{boyd_convex_optimization}, $V_t$ in \eqref{e:cost_compact} is convex.
	%	Now consider the objective function \eqref{e:cost_compact}.
	We substitute the stacked error vector \eqref{e:stacked error} in the objective function \eqref{e:cost_compact}.
	\begin{align*}
	& V_t  = \EE_{\XX_t} \left[ \norm{\calA e_t + \calB\control_{t:N}^e + \calD \wnoise_{t:N}}^2_{\calQ} + \norm{\control_{t:N}^e}^2_{\calR} \right] \\
	& = \EE_{\XX_t} \Bigl[ \norm{\calA e_t}^2_{\calQ} + \norm{\calD\wnoise_{t:N}}^2_{\calQ} + \norm{\control_{t:N}^e}^2_{\alpha} + 2( e_t \transp \calA \transp \calQ \calB  + \wnoise_{t:N} \transp \calD \transp \calQ \calB  ) \control_{t:N}^e  + 2 e_t\transp\calA\transp\calQ\calD\wnoise_{t:N} \Bigr].
	\end{align*}
	Let $\beta_t \Let \EE_{\XX_t}\left[ \norm{\calA e_t}^2_{\calQ} + \norm{\calD\wnoise_{t:N}}^2_{\calQ} + 2 e_t\transp\calA\transp\calQ\calD\wnoise_{t:N} \right].$ 
	Then $V_t = \EE_{\XX_t} \Bigl[\norm{\control_{t:N}^e}^2_{\alpha} + 2( e_t \transp \calA \transp \calQ \calB + \wnoise_{t:N} \transp \calD \transp \calQ \calB  ) \control_{t:N}^e \Bigr] + \beta_t$. We now substitute the stacked control vector \eqref{e:stacked_control_error} in $V_t$ to get
	\begin{align*}
V_t =	&  \EE_{\XX_t} \Bigl[\norm{(\calH -I)\control_{t:N}^r + \calG\offset_t + \calS \gain_t \ee(\tilde{\wnoise}_{t-1:N})}^2_{\alpha} + 2 \bigl( e_t \transp \calA \transp \calQ \calB  + \wnoise_{t:N} \transp \calD \transp \calQ \calB  \bigr) \bigl( (\calH -I)\control_{t:N}^r + \calG\offset_t + \calS \gain_t \ee(\tilde{\wnoise}_{t-1:N}) \bigr) \Bigr] + \beta_t.
	\end{align*}
	Let $\beta_t^{\prime} \Let \beta_t +  2 \EE_{\XX_t}\bigl[ e_t \transp \calA \transp \calQ \calB (\calH -I)\control_{t:N}^r]$. Since $\EE_{\XX_t}\left[ e_t \right] =  \stest_t -\st_t^r = e_t^{C}$, by removing mean zero terms we get the following equation:
	\begin{align}\label{e:objsolve}
	& V_t = \EE_{\XX_t} \Bigl[\norm{(\calH -I)\control_{t:N}^r + \calG\offset_t + \calS \gain_t \ee(\tilde{\wnoise}_{t-1:N})}^2_{\alpha}  + 2 \bigl( e_t \transp \calA \transp \calQ \calB + \wnoise_{t:N} \transp \calD \transp \calQ \calB \bigr) \calS\gain_t \ee(\tilde{\wnoise}_{t-1:N} ) \Bigr] + 2(e_t^C) \transp \calA \transp \calQ \calB \mu_{\calG}\offset_t +  \beta_t^{\prime}
	\end{align}
	Since $\tilde{\wnoise}_{t-1}$ is known at $t$, we simplify the term 
	\begin{align}
	& \EE_{\XX_t} \Bigl[\wnoise_{t:N} \transp \calD \transp \calQ \calB \calS \gain_t \ee(\tilde{\wnoise}_{t-1:N})\Bigr] \nonumber \\
	&= \EE_{\XX_t} \Biggl[ \wnoise_{t:N} \transp \calD \transp \calQ \calB \calS \bmat{\gain_t^{(:,t)} & \gain_t^{\prime}} \bmat{\ee(\tilde{\wnoise}_{t-1}) \\ \ee(\tilde{\wnoise}_{t:N-1})}\Biggr] = \trace \left( \calD \transp \calQ \calB\mu_{\calS}\gain_t^{\prime} \Sigma_{\ee\wnoise} \right) \label{e:SigmaPsi}.
	\end{align}
	We consider the term $\EE_{\XX_t} \Bigl[ e_ t \transp \calA \transp \calQ \calB \calS\gain_t \ee(\tilde{\wnoise}_{t-1:N})\Bigr]$ in \eqref{e:objsolve} as follows:
	\begin{align}\label{e:Sigmae}
	& \EE_{\XX_t} \Bigl[ e_ t \transp \calA \transp \calQ \calB \calS \gain_t \ee(\tilde{\wnoise}_{t-1:N})\Bigr] \notag \\
	& = \EE_{\XX_t} \Bigl[ (\st_ t - \stest_t) \transp \calA \transp \calQ \calB \calS \gain_t \ee(\tilde{\wnoise}_{t-1:N})\Bigr] + (\stest - \st_t^r)\transp \EE_{\XX_t} \Bigl[ \calA \transp \calQ \calB \calS \gain_t \ee(\tilde{\wnoise}_{t-1:N})\Bigr] \notag \\
	& = \EE_{\XX_t} \Bigl[ (e_t^D) \transp \calA \transp \calQ \calB \calS \gain_t \ee(\tilde{\wnoise}_{t-1:N})\Bigr] = \trace \Bigl( \calA \transp \calQ \calB \calS \gain_t^{\prime} \Sigma_{e\ee} \Bigr).
	\end{align}
	We substitute \eqref{e:SigmaPsi} and \eqref{e:Sigmae} in \eqref{e:objsolve} to get
	\begin{align}\label{e:reVt}
	& V_t = \EE_{\XX_t} \Bigl[\norm{(\calH -I)\control_{t:N}^r + \calG\offset_t + \calS \gain_t \ee(\tilde{\wnoise}_{t-1:N})}^2_{\alpha} \Bigr] \notag \\
	&   + 2 \trace \left( \calD \transp \calQ \calB\mu_{\calS}\gain_t^{\prime} \Sigma_{\ee\wnoise} \right) + 2\trace \Bigl( \calA \transp \calQ \calB \calS \gain_t^{\prime} \Sigma_{e\ee} \Bigr)  + 2(e_t^C) \transp \calA \transp \calQ \calB \mu_{\calG}\offset_t +  \beta_t^{\prime}
	\end{align}
	Let us define $c \Let (\control_{t:N}^r) \transp \EE_{\XX_t} \Bigl[ (\calH-I)\transp \alpha (\calH-I) \Bigr]\control_{t:N}^r$. We simplify the first term in the right hand side of \eqref{e:reVt} as follows:	
	\begin{align}\label{e:Vt}
	& \EE_{\XX_t} \Bigl[\norm{(\calH -I)\control_{t:N}^r + \calG\offset_t + \calS \gain_t \ee(\tilde{\wnoise}_{t-1:N})}^2_{\alpha} \Bigr] \notag \\
	& = \offset_t\transp \EE_{\XX_t} \Bigl[\calG\transp \alpha \calG \Bigr]\offset_t + \EE_{\XX_t} \Bigl[\norm{\calS\gain_t \ee(\tilde{\wnoise}_{t-1:N})}^2_{\alpha} + 2(\offset_t \transp \calG \transp \alpha) \calS\gain_t \ee(\tilde{\wnoise}_{t-1:N}) \Bigr] + 2\EE_{\XX_t} \left[(\control_{t:N}^r) \transp (\calH-I)\transp \alpha \calS \gain_t \psi(\tilde{\wnoise}_{t-1:N}) \right] \notag \\
	& \quad  + 2(\control_{t:N}^r) \transp \EE_{\XX_t} \left[(\calH-I)\transp \alpha \calG \right]\offset_t + \tilde{\beta}_t \notag \\
	& = \offset_t\transp \Sigma_{\calG} \offset_t + \EE_{\XX_t} \Bigl[\norm{\calS\gain_t \ee(\tilde{\wnoise}_{t-1:N})}^2_{\alpha} \Bigr] + 2\EE_{\XX_t} \Bigl[ \offset_t \transp \calG\transp\alpha\calS \gain_t \ee(\tilde{\wnoise}_{t-1:N}) \Bigr] + 2\EE_{\XX_t} \left[(\control_{t:N}^r) \transp (\calH-I)\transp \alpha \calS \gain_t \psi(\tilde{\wnoise}_{t-1:N}) \right] \notag \\
	& \quad + 2(\control_{t:N}^r) \transp \Sigma_{\calH\calG}\offset_t + c.
	\end{align}
	Let us consider the term $\EE_{\XX_t} \Bigl[\norm{\calS \gain_t \ee(\tilde{\wnoise}_{t-1:N})}^2_{\alpha} \Bigr]$ on the right hand side of \eqref{e:Vt}. In order to simplify offline computations, we perform the following manipulation:
	\begin{align}
	& \EE_{\XX_t} \Bigl[\norm{\calS\gain_t \ee(\tilde{\wnoise}_{t-1:N})}^2_{\alpha} \Bigr] \nonumber \\
	& = \EE_{\XX_t} \Biggl[\norm{\calS\bmat{\gain_t^{(:,t)} & \gain_t^{\prime}} \bmat{\ee(\tilde{\wnoise}_{t-1}) \\ \ee(\tilde{\wnoise}_{t:N-1})} }^2_{\alpha} \Biggr] \nonumber \\
	& = \trace \Biggl( \Sigma_{\calS} \gain_t^{(:,t)} \ee(\tilde{\wnoise}_{t-1}) \ee(\tilde{\wnoise}_{t-1}) \transp (\gain_t^{(:,t)})\transp + \Sigma_{\calS} \gain_t^{\prime} \EE_{\XX_t}\Bigl[\ee(\tilde{\wnoise}_{t:N-1})\ee(\tilde{\wnoise}_{t:N-1})\transp\Bigr] (\gain_t^{\prime})\transp    \Biggr) \nonumber \\
	& = \trace(\Sigma_{\calS} \gain_t^{(:,t)} \Pi_{\wnoise} (\gain_t^{(:,t)})\transp) + \trace (\Sigma_{\calS} \gain_t^{\prime} \Sigma_{\ee} (\gain_t^{\prime})\transp). \label{e:Sigmay} 
	\end{align}
	Let us consider the term $\EE_{\XX_t} \left[ \offset_t \transp \calG \transp \alpha \calS \gain_t \ee(\tilde{\wnoise}_{t-1:N})   \right]$ on the right hand side of \eqref{e:Vt} . By observing $\EE_{\XX_t} \left[ \ee( \tilde{\wnoise}_{t+i-1}) \right] = \zeros$ for each $i = 1, \ldots, N-1$, we get
	\begin{equation}\label{e:etay}
	\EE_{\XX_t} \left[ \offset_t \transp \calG\transp\alpha\calS \gain_t \ee(\tilde{\wnoise}_{t-1:N})   \right] = \offset_t \transp \Sigma_{\calG\calS}\gain_t^{(:,t)}\ee(\tilde{\wnoise}_{t-1}).
	\end{equation}
	Similar to \eqref{e:etay}, we get
	\begin{align}
	& \EE_{\XX_t} \left[(\control_{t:N}^r) \transp (\calH-I)\transp \alpha \calS \gain_t \psi(\tilde{\wnoise}_{t-1:N}) \right] = (\control_{t:N}^r) \transp \Sigma_{\calH\calS}\gain_t^{(:,t)}\ee(\tilde{\wnoise}_{t-1}). \label{e:etay ref}  
	\end{align}  
	Expression \eqref{e:obj out channel} follows by substituting \eqref{e:Vt}, \eqref{e:Sigmay},  \eqref{e:etay} and \eqref{e:etay ref} in \eqref{e:reVt}, and ignoring the term $\beta_t^{\prime} +c$, which is independent of the decision variables. 
	Therefore, the objective function in \eqref{e:cost_compact} is equivalent to \eqref{e:obj out channel} for the sake of optimization.  
\end{pf}
Before the proof of Theorem \ref{th:stability}, we need the following result related to mean square boundedness of $e_t^{D}$:
\begin{lem}\label{lem:msberror}
	Suppose that the dropout compensator is driven by the recursion \eqref{e:est} and let assumptions \ref{as:processnoise}-\ref{as:Lyapunov} hold, then there exists $\gamma^D > 0$ such that
	\begin{equation}\label{e:bound on dc}
	\EE_{\XX_{0}}\left[ \norm{e_{t}^{D}}^2 \right] \leq \gamma^D \quad \text{ for all } t .
	\end{equation}
\end{lem}
The proof of the above Lemma is along the lines of the proof of \cite[Lemma 9]{PDQ_intermittent}. Therefore, we omit the details for brevity.
\begin{pf}[Theorem \ref{th:stability}]
	Since $e_t^O = \st_t - r_t = \st_t - \stest_t + \stest_t - \st_t^r + \st_t^r - r_t = e_t^{D} + e_t^{C} + e_t^{G}$ and $\norm{e_t^O}^2 \leq 3 \left( \norm{e_t^{D}}^2 + \norm{e_t^{C}}^2 + \norm{e_t^{G}}^2 \right)$ by using Cauchy-Schwartz inequality. By taking conditional expectation on both sides we get \[ \EE_{\XX_0} \left[ \norm{e_t^O}^2 \right] \leq 3 \left( \EE_{\XX_0} \left[ \norm{e_t^{D}}^2 \right] + \EE_{\XX_0} \left[ \norm{e_t^{C}}^2 \right] + \norm{e_t^{G}}^2 \right).\] By \eqref{e:bound reference governor} we have
	\[
	\EE_{\XX_0} \left[ \norm{e_t^O}^2 \right]  \leq 3 \left( \EE_{\XX_0} \left[ \norm{e_t^{D}}^2 \right] + \EE_{\XX_0} \left[ \norm{e_t^{C}}^2 \right] + \gamma^{G} \right).
	\]
	Further, the Lemma \ref{lem:msberror} and the Lemma \ref{t:msbsingle} give us $ \EE_{\XX_0} \left[ \norm{e_t^O}^2 \right] \leq 3 \left( \gamma^{D} + \EE_{\XX_0} \left[ \norm{e_t^{C}}^2 \right] + \gamma^{G} \right)
	\leq 3 \left( \gamma^{D} + \gamma^{C} + \gamma^{G} \right) \teL \gamma $.
\end{pf}	
\end{document}